\documentclass[12pt]{amsart}
\usepackage{amsmath, amsthm, amssymb,cite}
\usepackage{fullpage}
\usepackage{color}
\usepackage[colorlinks=true]{hyperref}

\newtheorem{theorem}{Theorem}[section]
\newtheorem{lemma}{Lemma}

\newtheorem{definition}[lemma]{Definition}

\numberwithin{lemma}{section}

\allowdisplaybreaks

\numberwithin{equation}{section}

\makeatletter
\renewcommand{\subsection}{%
  \@startsection{subsection}{2}{\z@}%
  {1.5ex plus .5ex minus .2ex}%
  {.8ex plus .2ex}%
  {\normalfont\bfseries}%
}
\makeatother

\author{Jianxiang Liu}
\address{Department of Mathematics \& Statistics, McMaster University}
\curraddr{}
\email{liu2238@mcmaster.ca}

\begin{document}

\title{On the Existence of Nonclassical Solutions for Curvature Quotient Equations}

\begin{abstract}
For each $0\le l<k\le n$ with $k-l\ge3$, we construct a convex
Lipschitz viscosity solution of the prescribed curvature quotient
equation $\sigma_k(\kappa)/\sigma_l(\kappa)=f(x)$ with a smooth
positive right-hand side. The solution fails to have interior \(C^{1,\beta}\) regularity at the origin whenever \(1-\frac{2}{k-l}<\beta\le1\).
\end{abstract}

\maketitle
\tableofcontents

\section{Introduction}
Interior second-order estimates play a fundamental role in the
regularity theory of fully nonlinear elliptic equations.
Two important classes of equations in this context are Hessian
equations, involving the eigenvalues of $D^2u$, and curvature
equations, involving the principal curvatures of the graph of $u$.

Let $\Omega\subset\mathbb{R}^n$ be a domain. For $u\in C^2(\Omega)$,
write $\lambda(D^2u)$ for the eigenvalues of its Hessian and
$\kappa[u]$ for the principal curvatures of its graph $M_u=\{(x,u(x)):x\in\Omega\}
    \subset\mathbb R^{n+1}$. Let $\sigma_j$ denote the $j$th elementary
symmetric polynomial, with $\sigma_0=1$. For integers $0\le l<k\le n$,
the Hessian quotient equation takes the form
\begin{align*}
    \left(\frac{\sigma_k(\lambda(D^2u))}
                    {\sigma_l(\lambda(D^2u))}\right)^{1/(k-l)}
    = g.
\end{align*}
In this paper we study the corresponding prescribed curvature
quotient equation
\begin{align*}
    G_{k,l}[u]
    := \left(\frac{\sigma_k(\kappa[u])}
                  {\sigma_l(\kappa[u])}\right)^{1/(k-l)}
    = g,
\end{align*}
where $g$ is a prescribed positive function.

Let $\Gamma_k=\{\lambda\in\mathbb{R}^n:
\sigma_j(\lambda)>0,\ 1\le j\le k\}$. A solution $u\in C^2(\Omega)$
is called $k$-admissible for the Hessian or curvature equation if
$\lambda(D^2u)\in\Gamma_k$ or $\kappa[u]\in\Gamma_k$, respectively,
throughout $\Omega$. In particular, $\Gamma_n=\Gamma_+=(0,\infty)^n$.

The interior regularity of these equations has been extensively investigated under various convexity and admissibility assumptions. We first recall some developments concerning interior second-order estimates for Hessian equations and Hessian quotient equations.

Interior Hessian estimates for the two-dimensional Monge--Amp\`ere equation were established by Heinz~\cite{Heinz1959}. For the quadratic Hessian equation $\sigma_2(D^2u)=f$, significant progress has been made in dimensions three and four~\cite{WarrenYuan2009,warrenYuan2010,Qiu2024hessian,Zhou2024,ShankarYuan2025,Fan2026D4}. In higher dimensions, interior estimates under various convexity or integrability assumptions were obtained in~\cite{GuanQiu2019,McGonagleSongYuan2019,ShankarYuan2020,mooney2021,mooney2025,chenjianz2026}. More recently, Li and Wu~\cite{LiWu2026Sigma2} established interior Hessian estimates for admissible solutions of $\sigma_2(D^2u)=1$ in arbitrary dimensions without additional convexity assumptions. Related regularity results for viscosity solutions were obtained in~\cite{shankaryuan2021,mooney2021,chenjiantuzhou2026,chenzhouzhu2026,zhouzhu2026}.

Interior Hessian estimates for Hessian quotient equations have also attracted considerable attention. Early results for $\sigma_3(D^2u)/\sigma_1(D^2u)$ in dimensions three and four were obtained by exploiting its special Lagrangian structure~\cite{ChenYuanWarren2009,WangYuan2014}, while Zhou~\cite{Zhou2024} treated variable right-hand sides using a twisted Lagrangian structure. For the quotients $\sigma_n(D^2u)/\sigma_{n-i}(D^2u)$, $i=1,2$, Lu~\cite{Lu2023,Lu2025} established interior estimates for convex solutions, subsequently extended to gradient-dependent right-hand sides under an additional concavity assumption by Fung~\cite{Fung2026}. More recently, significant progress has been made for general Hessian quotient equations with $k-l=1,2$, under various convexity, semiconvexity, and structural assumptions~\cite{lutsai2026note,jiaosui2026,lusroka2026,jiliang2026,LiWu2026Quotient,DongZhang2026}. For constant right-hand sides, Wu~\cite{Wu2026}
obtained interior Hessian estimates for $k$-admissible
solutions when $2k>n$, while Qiu and Yan~\cite{QiuYan2026}
established such estimates in arbitrary dimensions.
Li and Wu~\cite{liwu2026generalquotient}
subsequently provided an alternative proof based on
a compactness argument.

For prescribed curvature equations, Urbas~\cite{Urbas2000Curvature}
established interior curvature estimates under additional
integrability assumptions on the mean curvature.
For the scalar curvature equation, Guan and Qiu~\cite{GuanQiu2019}
obtained interior curvature estimates for convex graphs, while
Qiu~\cite{Qiu2024Curvature} treated admissible graphs in dimension
three. Further developments include the estimates for convex
graphs with positive Lipschitz right-hand sides
by Chen, Jian, and Zhou~\cite{chenjianz2026}, and the
four-dimensional results of Fan and Shankar~\cite{FanShankar2026}.
More recently, Qiu and Yan~\cite{QiuYan2026Curvature}
established interior curvature estimates for admissible
solutions of $\sigma_2(\kappa[u])=1$ in arbitrary dimensions
$n\ge3$, without additional convexity assumptions.

For curvature quotient equations, Sheng, Urbas, and Wang
\cite{ShengUrbasWang2004} established Pogorelov-type
interior curvature estimates under affine boundary conditions.
Qiu and Zhou~\cite{QiuZhou2024} obtained interior curvature
estimates for special Lagrangian curvature equations,
which also cover certain curvature quotient equations
in low dimensions.
Liu~\cite{Liu2026} established interior curvature
estimates for convex graphical solutions of
$\sigma_n(\kappa[u])/\sigma_{n-2}(\kappa[u])=f(X)$
in arbitrary dimensions.
Subsequently, Wang~\cite{Wang2026CurvatureQuotient} obtained interior curvature
estimates for a class of elliptic Weingarten equations
on strictly convex graphs, including certain curvature
quotient equations.
More recently, Han and Ma~\cite{HanMa2026} established
interior curvature estimates for admissible graphical solutions of
$\sigma_k(\kappa[u])/\sigma_{k-2}(\kappa[u])=1$,
$3\le k<n$, without additional convexity or
semiconvexity assumptions.

Despite the substantial progress described above, interior estimates
do not always hold for Hessian equations and prescribed curvature
equations. The classical singular examples of
Pogorelov~\cite{Pogorelov1978} show that such estimates fail in general
for $\det D^2u=\sigma_3(D^2u)=1$ in dimension three.
Urbas~\cite{Urbas1990} extended this phenomenon to higher-order
Hessian equations $\sigma_k(D^2u)=f$, $n\geq k\geq 3$.
Following the approach of Urbas, Lu~\cite[Theorem 1.2]{Lu2025}
constructed counterexamples showing that interior estimates
fail in general for Hessian quotient equations
$\sigma_k(D^2u)/\sigma_l(D^2u)=f$ when $1\leq l<k\leq n$
and $k-l\geq 3$, even for convex solutions.
Furthermore, Fung~\cite[Remark 1.4]{Fung2026}, using the idea
of~\cite[Theorem 1.2]{GuanRenWang2015}, showed that such estimates
may fail even when $k-l=1,2$ if the right-hand side depends on
$Du$ without appropriate structural assumptions.
These results reveal the critical role of the degree gap $k-l$
in the interior regularity of Hessian quotient equations.

For prescribed curvature quotient equations,
Han and Ma~\cite{HanMa2026} pointed out that the restriction
$k-l\leq 2$ is optimal at least in the case $l=0$,
in view of the counterexamples constructed by
Urbas~\cite{Urbas1990}. 
In this paper, we show that interior curvature estimates fail in general for all curvature quotient equations with $1\leq l<k\leq n$ and $k-l\geq 3$, even for convex graphical solutions with smooth positive right-hand
sides depending only on $x$ by constructing nonclassical viscosity solutions.

\begin{theorem}\label{thm:main}
Let $n\ge k\ge 3$ and $0\le l\le k-3$ be integers.
Then there exist $\varepsilon\in(0,1)$, a smooth positive function
$g$ on $B_1\subset\mathbb{R}^n$, and a convex viscosity solution
$u$ of
\begin{align}\label{eq:k,lquotient}
    G_{k,l}[u]=g(x) \qquad\text{in } B_\varepsilon
\end{align}
such that
\begin{align*}
    u\in C^{0,1}(\overline{B_\varepsilon}),\qquad
    u\notin C^{1,\beta}(B_{\varepsilon/2})
    \quad\text{for every } 1-\frac{2}{k-l}<\beta\le 1.
\end{align*}
\end{theorem}

The proof follows the strategy of Urbas~\cite[Theorem~2]{Urbas1990}.
In Section~\ref{sec:preliminaries}, we extend the graph curvature operator
and construct the barriers used in the approximating Dirichlet problems.
In Section~\ref{sec:proof1.1}, we solve these Dirichlet problems, use the barriers
to obtain uniform Lipschitz bounds, and pass to a viscosity limit.
The growth of the limit at the origin gives the asserted failure of
regularity. Constant extension in the remaining variables yields the case $k<n$.

\section{Preliminaries}\label{sec:preliminaries}

\subsection{Viscosity solutions and the graph curvature operator}

We use the sign convention of~\cite[Definition~2.3]{CaffarelliCabre1995}
in the following definition.

\begin{definition}\label{def:viscosity}
Let $F\colon \Omega\times\mathbb R\times\mathbb R^n\times\mathbb S^n
\to\mathbb R$ be continuous and degenerate elliptic, where
$\Omega\subset\mathbb R^n$ is a domain and $\mathbb S^n$ denotes
the space of real symmetric $n\times n$ matrices. A function $u\in C(\Omega)$ is a \emph{viscosity subsolution}
(respectively, \emph{viscosity supersolution}) of
\begin{align*}
    F(x,u,Du,D^2u)=0 \qquad\text{in }\Omega
\end{align*}
if, for every $x_0\in\Omega$ and $\varphi\in C^2(\Omega)$ such that
$u-\varphi$ has a local maximum (respectively, local minimum)
at $x_0$,
\begin{align*}
    F(x_0,u(x_0),D\varphi(x_0),D^2\varphi(x_0))
    \ge0
    \qquad(\text{respectively, }\le0).
\end{align*}
A \emph{viscosity solution} is both a viscosity subsolution and
a viscosity supersolution.
\end{definition}

We now write the curvature quotient equation $G_{k,l}[u]=g(x)$
in terms of $Du$ and $D^2u$, where $k,l$ take the same values as in \ref{thm:main}.
By \cite{CaffarelliNirenbergSpruck1988}, the principal curvatures
$\kappa[u]$ of the graph of $u\in C^2(\Omega)$ are the eigenvalues
of the symmetric matrix
\begin{align}\label{eq:graphmatrix}
    a[u]=\frac{1}{w}\gamma D^2u\,\gamma,\qquad
    w=\sqrt{1+|Du|^2},\qquad
    \gamma=\mathsf g^{-1/2}
    =I-\frac{Du\otimes Du}{w(1+w)},
\end{align}
where $\mathsf g=I+Du\otimes Du$ is the induced metric.
For the quotient $G_{k,l}$, denote
\begin{align*}
    G_{k,l}[u]=f(\kappa[u]),\qquad
    f(\kappa)=
    \left(\frac{\sigma_k(\kappa)}{\sigma_l(\kappa)}\right)^{1/(k-l)},
    \qquad \kappa\in\Gamma_+.
\end{align*}

In particular, for $k=n$, writing
$\kappa^{-1}=(\kappa_1^{-1},\ldots,\kappa_n^{-1})$, we have
\begin{align}\label{eq:reciprocal}
    \frac{\sigma_l(\kappa)}{\sigma_n(\kappa)}
    =\sigma_{n-l}(\kappa^{-1}).
\end{align}

\subsection{Continuous extension and strict comparison}\label{sec:lmconcheck}
For the remainder of this section, we restrict to the case $k=n$.
Extend $f$ to $\mathbb R^n$ by
\begin{align}\label{eq:fhat}
    \widehat f(\kappa)=
    \begin{cases}
        f(\kappa), & \kappa\in\Gamma_+,\\
        0,         & \kappa\notin\Gamma_+.
    \end{cases}
\end{align}
For $0<\kappa_1\le\cdots\le\kappa_n$, \eqref{eq:reciprocal} gives
\begin{align*}
    0<f(\kappa)^{n-l}
    =\frac{1}{\sigma_{n-l}(\kappa^{-1})}
    \le \kappa_1\kappa_n^{n-l-1}.
\end{align*}
At each finite point of $\partial\Gamma_+$, the right-hand side tends to
zero since $\kappa_1\to0$ and $\kappa_n$ remains bounded,
thus $\widehat f\in C(\mathbb R^n)$.
We still denote $G_{n,l}[u]=\widehat f(\kappa[u])$.

Let $\Omega\Subset B_1$ and let $g$ be smooth and positive on $B_1$.
Write $a(p,R)$ for the matrix in \eqref{eq:graphmatrix} with
$(Du,D^2u)=(p,R)$, and $\kappa(a(p,R))$ for its ordered eigenvalues.
Define
\begin{align}\label{eq:Q}
    Q(x,z,p,R)=\widehat f\bigl(\kappa(a(p,R))\bigr)-g(x).
\end{align}
Then $Q$ is continuous. For $t\ge0$ and $x,y\in\overline\Omega$,
\begin{align*}
    Q(x,z,p,R)-Q(x,z+t,p,R) &= 0,\\
    |Q(x,z,p,R)-Q(y,z,p,R)| &= |g(x)-g(y)|
    \le \omega_g(|x-y|),
\end{align*}
where $\omega_g$ is a modulus of continuity of $g$ on
$\overline\Omega$. Hence $Q$ and its constant translates satisfy
conditions (2.4) and (2.5) of \cite{Trudinger1990} for
comparison principle~\cite[Lemma~2.1]{Trudinger1990}. It is also easy to check $Q$ is degenerate elliptic, which is the condition (10) in \cite[Proposition I.3]{Lions1983} (we use $-Q$ here).

\subsection{Barrier estimates}

We construct $V_\sigma$ as a subsolution for the approximating Dirichlet problem and control the gradient, and $\psi_\sigma$ to provide a lower bound for the growth of the limit at the origin.
Fix a smooth positive function $g$ on $B_1$. Set
\[
 \alpha=2-\frac2{n-l},\qquad x'=(x_2,\ldots,x_n),
 \qquad s=\sigma+|x'|^2.
\]
Since $n-l\ge3$, we have $4/3\le\alpha<2$. For $A>0$ and
$0<\sigma\le\varepsilon^2$, define
\begin{equation}\label{eq:barriers}
 \varphi_\sigma=2As^{\alpha/2},\qquad
 \psi_\sigma=A(1+x_1^2)s^{\alpha/2},\qquad
  V_\sigma=\varphi_\sigma+\frac{|x|^2-\varepsilon^2}{\varepsilon}.
\end{equation}
Choose
\begin{equation}\label{eq:Achoice}
 A>9\cdot2^{3/2}\binom n{n-l}^{1/(n-l)}(\sup_{\overline{B_{1/2}}} g+1),
\end{equation}
and then choose
\begin{equation}\label{eq:epschoice}
 0<\varepsilon\le\min\left\{
 \frac18,\ (8A)^{-1/(\alpha-1)},\
 \frac{2}{10^{3/2}\binom n{n-l}^{1/(n-l)}(\sup_{\overline{B_{1/2}}} g+1)}
 \right\}.
\end{equation}
These choices are independent of $\sigma$. We then have the following two lemmas.

\begin{lemma}\label{lem:Vbarrier}
For every $0<\sigma\le\varepsilon^2$, $V_\sigma$ is smooth and strictly
convex, $\|DV_\sigma\|_{L^\infty(B_\varepsilon)}\le3$, and
\[
 G_{n,l}[ V_\sigma]\ge \sup_{\overline{B_{1/2}}} g+1\quad\text{in }B_\varepsilon,\qquad
  V_\sigma=\varphi_\sigma\quad\text{on }\partial B_\varepsilon.
\]
\end{lemma}

\begin{proof}
The boundary equality follows from \eqref{eq:barriers}. The gradient is
\[
 D\varphi_\sigma=(0,2A\alpha x_2s^{\alpha/2-1},\ldots,
 2A\alpha x_ns^{\alpha/2-1}).
\]
For $r=|x'|>0$, since $\alpha/2-1<0$,
\[
 |D\varphi_\sigma|=2A\alpha r(\sigma+r^2)^{\alpha/2-1}
 \le2A\alpha r^{\alpha-1}\le2A\alpha\varepsilon^{\alpha-1}.
\]
At $r=0$, the gradient formula gives $D\varphi_\sigma=0$ because $s=\sigma>0$.
Consequently,
\begin{align}
 \|D\varphi_\sigma\|_{L^\infty(B_\varepsilon)}
 &\le2A\alpha\varepsilon^{\alpha-1}\le1,\label{eq:gradphi}\\
 \|D V_\sigma\|_{L^\infty(B_\varepsilon)}
 &\le2+2A\alpha\varepsilon^{\alpha-1}\le3,\label{eq:gradV}
\end{align}
where the last inequality for both formulae comes from \eqref{eq:Achoice} and \eqref{eq:epschoice}.

For $i,j\ge2$,
\[
 D_{ij}\varphi_\sigma=2A\alpha s^{\alpha/2-1}\delta_{ij}
 +2A\alpha(\alpha-2)s^{\alpha/2-2}x_ix_j,
\]
and $D_{1j}\varphi_\sigma=D_{j1}\varphi_\sigma=0$ for all $j$. For $\xi=(\xi_1,\xi')$,
Cauchy--Schwarz and $\alpha-2<0$ yield
\begin{align*}
 D^2\varphi_\sigma[\xi,\xi]
 &=2A\alpha s^{\alpha/2-2}
 \bigl[s|\xi'|^2+(\alpha-2)(x'\cdot\xi')^2\bigr]\\
 &\ge2A\alpha s^{\alpha/2-2}
 \bigl[s+(\alpha-2)|x'|^2\bigr]|\xi'|^2\\
 &=2A\alpha s^{\alpha/2-2}
 \bigl[\sigma+(\alpha-1)|x'|^2\bigr]|\xi'|^2\ge0.
\end{align*}
Thus
\begin{equation}\label{eq:hessianV}
 D^2 V_\sigma=D^2\varphi_\sigma+\frac2\varepsilon I\ge\frac2\varepsilon I>0.
\end{equation}
Apply \eqref{eq:graphmatrix} with $u=V_\sigma$.
Since $\mathsf g\le w^2I$, we have $\gamma^2\ge w^{-2}I$.
For every unit vector $\xi$,
\begin{align*}
 \xi^Ta[ V_\sigma]\xi
 &=\frac1w(\gamma\xi)^TD^2 V_\sigma(\gamma\xi)\\
 &\ge\frac{\kappa_{\min}(D^2 V_\sigma)}w|\gamma\xi|^2
 \ge\frac{\kappa_{\min}(D^2 V_\sigma)}{w^3}.
\end{align*}
By \eqref{eq:gradV} and \eqref{eq:hessianV},
\begin{equation}\label{eq:curvatureV}
 \kappa_{\min}(a[ V_\sigma])\ge\frac2{10^{3/2}\varepsilon}.
\end{equation}
Using \eqref{eq:reciprocal}, we obtain
\begin{align}
 G_{n,l}[ V_\sigma]&=\bigl[\sigma_{n-l}(\kappa[ V_\sigma]^{-1})\bigr]^{-1/(n-l)}\notag\\
 &\ge\binom n{n-l}^{-1/(n-l)}\kappa_{\min}(a[ V_\sigma])\notag\\
 &\ge\frac2{10^{3/2}\binom n{n-l}^{1/(n-l)}\varepsilon}.
 \label{eq:FV}
\end{align}
The choice \eqref{eq:epschoice} gives the asserted lower bound.
\end{proof}

\begin{lemma}\label{lem:psibarrier}
For every $0<\sigma\le\varepsilon^2$,
\[
 G_{n,l}[\psi_\sigma]\ge\sup_{\overline{B_{1/2}}} g+1\quad\text{in }B_\varepsilon,
 \qquad \psi_\sigma\le\varphi_\sigma\quad\text{on }\partial B_\varepsilon.
\]
\end{lemma}

\begin{proof}
The first derivatives are
\[
 D_1\psi_\sigma=2Ax_1s^{\alpha/2},\qquad
 D_j\psi_\sigma=A\alpha(1+x_1^2)x_js^{\alpha/2-1}\quad(j\ge2).
\]
Since $0<\sigma\le\varepsilon^2$ and $\varepsilon\le1$,
\begin{align}
 |D\psi_\sigma|&\le2A|x_1|s^{\alpha/2}
 +A\alpha(1+x_1^2)|x'|s^{\alpha/2-1}\notag\\
 &\le2A\varepsilon(2\varepsilon^2)^{\alpha/2}
 +A\alpha(1+\varepsilon^2)\varepsilon^{\alpha-1}
 \le8A\varepsilon^{\alpha-1}.\label{eq:gradpsi}
\end{align}
The second derivatives are
\begin{align}
 D_{11}\psi_\sigma&=2As^{\alpha/2},\notag\\
 D_{1j}\psi_\sigma&=2A\alpha x_1x_js^{\alpha/2-1}\quad(j\ge2),\notag\\
 D_{ij}\psi_\sigma&=A\alpha(1+x_1^2)s^{\alpha/2-1}\delta_{ij}
 +A\alpha(\alpha-2)(1+x_1^2)s^{\alpha/2-2}x_ix_j
 \quad(i,j\ge2).\label{eq:hessianpsi}
\end{align}
For $\xi=(\xi_1,\xi')$, since $|x'|^2\le s$ and $\alpha-2<0$,
\begin{align*}
 A^{-1}D^2\psi_\sigma[\xi,\xi]
 &\ge2s^{\alpha/2}\xi_1^2
 -4\alpha|x_1|s^{\alpha/2-1}|x'||\xi_1||\xi'|\\
 &\hspace{1cm}+\alpha(\alpha-1)(1+x_1^2)s^{\alpha/2-1}|\xi'|^2.
\end{align*}
Young's inequality gives
\[
 4\alpha|x_1|s^{\alpha/2-1}|x'||\xi_1||\xi'|
 \le s^{\alpha/2}\xi_1^2+4\alpha^2x_1^2s^{\alpha/2-1}|\xi'|^2.
\]
It follows that
\[
 A^{-1}D^2\psi_\sigma[\xi,\xi]\ge s^{\alpha/2}\xi_1^2
 +[\alpha(\alpha-1)-4\alpha^2\varepsilon^2]s^{\alpha/2-1}|\xi'|^2.
\]
Because $4/3\le\alpha<2$ and $\varepsilon\le1/8$,
\[
 \alpha(\alpha-1)-4\alpha^2\varepsilon^2
 \ge\frac49-16\varepsilon^2\ge\frac7{36}>\frac19.
\]
Therefore
\begin{align}
 D^2\psi_\sigma&\ge\frac A9
 \operatorname{diag}(s^{\alpha/2},s^{\alpha/2-1},\ldots,s^{\alpha/2-1})>0,
 \label{eq:Hlower}\\
 (D^2\psi_\sigma)^{-1}&\le\frac9A
 \operatorname{diag}(s^{-\alpha/2},s^{1-\alpha/2},\ldots,s^{1-\alpha/2}).
 \label{eq:Hinverse}
\end{align}
Applying $\sigma_{n-l}$ to \eqref{eq:Hinverse} yields
\begin{align}
 \sigma_{n-l}\bigl((D^2\psi_\sigma)^{-1}\bigr)
 &\le\left(\frac9A\right)^{n-l}
 \left[\binom{n-1}{n-l-1}s^{-\alpha/2+(n-l-1)(1-\alpha/2)}
 +\binom{n-1}{n-l}s^{(n-l)(1-\alpha/2)}\right]\notag\\
 &=\left(\frac9A\right)^{n-l}
 \left[\binom{n-1}{n-l-1}+\binom{n-1}{n-l}s\right]\notag\\
 &\le\left(\frac9A\right)^{n-l}\binom n{n-l}.
 \label{eq:inversebound}
\end{align}
Here $\alpha/2=1-1/(n-l)$ and $s\le2\varepsilon^2\le1$; the convention
$\binom{n-1}{n}=0$ includes the case $l=0$.

To pass from the Hessian to principal curvatures, put $H=D^2\psi_\sigma$ and let
\[
 0<\mu_1\le\cdots\le\mu_n,\qquad
 0<\kappa_1\le\cdots\le\kappa_n
\]
be the eigenvalues of $H$ and $a[\psi_\sigma]=w^{-1}\gamma H\gamma$, respectively.
For $\xi\ne0$, put $\eta=\gamma\xi$. Since
$\gamma^{-2}=\mathsf g\le w^2I$,
\begin{equation}\label{eq:rayleigh}
 \frac{\xi^Ta[\psi_\sigma]\xi}{|\xi|^2}
 =\frac{\eta^TH\eta}{w\eta^T\mathsf g\eta}
 \ge\frac1{w^3}\frac{\eta^TH\eta}{|\eta|^2}.
\end{equation}
For increasingly ordered eigenvalues, the min--max principle gives
\[
 \kappa_i=\min_{\substack{E\subset\mathbb{R}^n\\\dim E=i}}
 \max_{0\ne\xi\in E}\frac{\xi^Ta[\psi_\sigma]\xi}{|\xi|^2}.
\]
Because $\gamma$ is invertible, $\gamma E$ runs through all
$i$-dimensional subspaces as $E$ does. Hence \eqref{eq:rayleigh} implies
\begin{equation}\label{eq:eigencomparison}
 \kappa_i\ge\frac1{w^3}
 \min_{\substack{L\subset\mathbb{R}^n\\\dim L=i}}
 \max_{0\ne\eta\in L}\frac{\eta^TH\eta}{|\eta|^2}
 =\frac{\mu_i}{w^3},\qquad 1\le i\le n.
\end{equation}
Consequently,
\begin{align*}
 \sigma_{n-l}(\kappa^{-1})
 &=\sum_{i_1<\cdots<i_{n-l}}
 \frac1{\kappa_{i_1}\cdots\kappa_{i_{n-l}}}\\
 &\le w^{3(n-l)}\sum_{i_1<\cdots<i_{n-l}}
 \frac1{\mu_{i_1}\cdots\mu_{i_{n-l}}}\\
 &=w^{3(n-l)}\sigma_{n-l}(H^{-1})\\
 &\le w^{3(n-l)}\left(\frac9A\right)^{n-l}\binom n{n-l}.
\end{align*}
By \eqref{eq:reciprocal} and \eqref{eq:gradpsi},
\begin{align}
 G_{n,l}[\psi_\sigma]&=\sigma_{n-l}(\kappa^{-1})^{-1/(n-l)}\notag\\
 &\ge\frac{A}{9w^3\binom n{n-l}^{1/(n-l)}}\notag\\
 &\ge\frac{A}{9(1+64A^2\varepsilon^{2\alpha-2})^{3/2}
 \binom n{n-l}^{1/(n-l)}}.\label{eq:Fpsi}
\end{align}
By \eqref{eq:Achoice} and \eqref{eq:epschoice},
\begin{equation}\label{eq:strictpsi}
 G_{n,l}[\psi_\sigma]\ge\frac{A}{9\cdot2^{3/2}\binom n{n-l}^{1/(n-l)}}
 >\sup_{\overline{B_{1/2}}}g+1.
\end{equation}
On $\partial B_\varepsilon$, since $x_1^2\le\varepsilon^2\le1$,
\[
 \psi_\sigma=A(1+x_1^2)s^{\alpha/2}\le2As^{\alpha/2}=\varphi_\sigma.
\]
\end{proof}

\section{Proof of Theorem \ref{thm:main}}\label{sec:proof1.1}
We first prove the case $k=n$. Choose $A$ and then $\varepsilon$
according to \eqref{eq:Achoice} and \eqref{eq:epschoice}, respectively.

\subsection*{Step 1. Existence of the approximating solutions}
For each $0<\sigma\le\varepsilon^2$, consider the Dirichlet problem
\begin{equation}\label{eq:dirichlet}
 \begin{cases}
 G_{n,l}[ u_\sigma]=g(x)&\text{in }B_\varepsilon,\\
  u_\sigma=\varphi_\sigma&\text{on }\partial B_\varepsilon.
 \end{cases}
\end{equation}
By \cite[Remark~4.1]{ShengUrbasWang2004}, the operator $G_{n,l}$
satisfies the structural hypotheses of Theorem~4.2 therein.
The ball $B_\varepsilon$ is smooth and uniformly convex. For each
$0<\sigma\le\varepsilon^2$, we have
$\varphi_\sigma,g\in C^\infty(\overline{B_\varepsilon})$,
with $g>0$ and $g_z=0$.
Thus, to obtain a classical convex solution $u_\sigma$ of
\eqref{eq:dirichlet} by \cite[Theorem~4.2]{ShengUrbasWang2004},
it suffices to exhibit a convex subsolution with boundary values
$\varphi_\sigma$. Lemma~\ref{lem:Vbarrier} provides such a subsolution $V_\sigma$.
Thus for each $0<\sigma\le\varepsilon^2$ there exists a unique
convex solution
\[
  u_\sigma\in C^{3,\theta}(B_\varepsilon)\cap C^{0,1}(\overline{B_\varepsilon}),
 \qquad 0<\theta<1,
\]
of \eqref{eq:dirichlet}. Since $Q[V_\sigma]\ge 1$, $Q[u_\sigma]=0$ in $B_\varepsilon$, and
$V_\sigma=u_\sigma$ on $\partial B_\varepsilon$, the strict comparison
principle \cite[Lemma~2.1]{Trudinger1990} gives
\begin{align}\label{eq:Vbelowu}
    V_\sigma\le u_\sigma \quad\text{in } B_\varepsilon.
\end{align}

\subsection*{Step 2. Comparison bounds and uniform Lipschitz estimates}
Since $\varphi_\sigma$ is independent of $x_1$, its curvature matrix
has a zero eigenvalue and $G_{n,l}[\varphi_\sigma]=0$ by \eqref{eq:fhat}.
Let $\mu_0=\min_{\overline{B_{1/2}}}g>0$ and $Q_1=Q+\mu_0/2$. Then
\[
 Q_1[ u_\sigma]=\frac{\mu_0}2,\qquad
 Q_1[\varphi_\sigma]=-g+\frac{\mu_0}2\le-\frac{\mu_0}2,
 \qquad  u_\sigma=\varphi_\sigma\quad\text{on }\partial B_\varepsilon.
\]
Strict comparison principle~\cite[Lemma~2.1]{Trudinger1990} gives
$u_\sigma\le\varphi_\sigma$. For the lower bound,
Lemma~\ref{lem:psibarrier} gives $Q[\psi_\sigma]\ge1$ in $B_\varepsilon$
and $\psi_\sigma\le\varphi_\sigma=u_\sigma$ on $\partial B_\varepsilon$.
Again, applying strict comparison principle to $Q$ yields
\begin{equation}\label{eq:sandwich}
 0\le\psi_\sigma\le u_\sigma\le\varphi_\sigma\quad\text{in }B_\varepsilon,
 \qquad 0<\sigma\le\varepsilon^2.
\end{equation}
In particular,
\begin{equation}\label{eq:C0}
 \| u_\sigma\|_{L^\infty(B_\varepsilon)}
 \le2A(2\varepsilon^2)^{\alpha/2}=2^{1+\alpha/2}A\varepsilon^\alpha.
\end{equation}
For the uniform gradient bound, use \eqref{eq:Vbelowu} and the common
boundary values. Given $x\in B_\varepsilon$ and a unit vector $e$, choose
$y=x+te\in\partial B_\varepsilon$ with $t>0$. By convexity,
\begin{align*}
 tD u_\sigma(x)\cdot e&\le u_\sigma(y)- u_\sigma(x)\\
 &\le V_\sigma(y)- V_\sigma(x)\le\|D V_\sigma\|_{L^\infty(B_\varepsilon)}t.
\end{align*}
As $e$ is arbitrary, \eqref{eq:gradV} implies
\begin{equation}\label{eq:C1}
 \|D u_\sigma\|_{L^\infty(B_\varepsilon)}
 \le2+2A\alpha\varepsilon^{\alpha-1}\le3.
\end{equation}
This Lipschitz bound extends to the closed ball by continuity.

\subsection*{Step 3. The viscosity limit and failure of regularity}
By \eqref{eq:C0}, \eqref{eq:C1}, and the Arzel\`a--Ascoli theorem,
there is a sequence $\sigma_j\downarrow0$ such that
\begin{equation}\label{eq:uniformlimit}
 u_{\sigma_j}\longrightarrow u_0\quad\text{uniformly on }\overline{B_\varepsilon},
 \qquad u_0\in C^{0,1}(\overline{B_\varepsilon}).
\end{equation}
The limit $u_0$ is convex, with the same uniform Lipschitz bound.

The operator $\mathcal F=-Q$ from Section~\ref{sec:preliminaries} is independent of $\sigma$.
Each classical solution $u_{\sigma_j}$ is a viscosity solution of
$\mathcal F=0$. \cite[Proposition I.3]{Lions1983} and \eqref{eq:uniformlimit}
therefore show that $u_0$ is a convex viscosity solution of
\begin{equation}\label{eq:limitproblem}
 \begin{cases}
 G_{n,l}[u_0]=g(x)&\text{in }B_\varepsilon,\\
 u_0=\varphi_0&\text{on }\partial B_\varepsilon,
 \end{cases}
 \qquad \varphi_0=2A|x'|^\alpha.
\end{equation}
The boundary equality follows from convergence on the closed ball and
$\varphi_\sigma\to\varphi_0$ uniformly. Similarly,
\[
 \psi_\sigma\longrightarrow\psi_0=A(1+x_1^2)|x'|^\alpha
 \quad\text{uniformly on }\overline{B_\varepsilon}.
\]
For completeness, since $0<\alpha/2<1$,
\[
 0\le(\sigma+|x'|^2)^{\alpha/2}-|x'|^\alpha\le\sigma^{\alpha/2},
\]
which verifies both barrier convergences directly.
Passing to the limit in \eqref{eq:sandwich} gives
\begin{equation}\label{eq:limitsandwich}
 0\le\psi_0\le u_0\le\varphi_0\quad\text{in }B_\varepsilon.
\end{equation}
Both barriers vanish on $B_\varepsilon\cap\{x_2=\cdots=x_n=0\}$.
In particular, $u_0(0)=0$ and $u_0\ge0$. Along $x=te_2$,
\begin{equation}\label{eq:axis}
 A|t|^{2-2/(n-l)}\le u_0(te_2)\le2A|t|^{2-2/(n-l)},
 \qquad |t|<\varepsilon.
\end{equation}
Suppose that $u_0\in C^{1,\beta}(B_{\varepsilon/2})$ for some
$\beta>1-2/(n-l)$. The origin is a minimum point, so $Du_0(0)=0$.
Thus, for sufficiently small $|t|$,
\[
 A|t|^{2-2/(n-l)}\le u_0(te_2)\le C|t|^{1+\beta}.
\]
Hence
\[
 0<A\le C|t|^{\beta-1+2/(n-l)}\longrightarrow0\qquad(t\to0),
\]
a contradiction to \eqref{eq:Achoice}. Therefore $u_0\notin C^{1,\beta}(B_{\varepsilon/2})$ for every
$\beta>1-2/(n-l)$. This proves the case $k=n$.

For $k<n$, the preceding construction
in dimension $k$ gives classical convex solutions $u_{\sigma_j}$ of
$G_{k,l}[u_{\sigma_j}]=g_0$ in $B_\varepsilon^k$ converging uniformly
on $\overline{B_\varepsilon^k}$ to $u_0$.
For $(x,y)\in\mathbb R^k\times\mathbb R^{n-k}$, set
$\widetilde u_{\sigma_j}(x,y)=u_{\sigma_j}(x)$,
$\widetilde u_0(x,y)=u_0(x)$, and $\widetilde g(x,y)=g_0(x)$,
so that $\widetilde g$ is smooth and positive on $B_1^n$.
By \eqref{eq:graphmatrix},
\begin{align*}
    a[\widetilde u_{\sigma_j}]
    =\begin{pmatrix} a[u_{\sigma_j}] & 0 \\ 0 & 0 \end{pmatrix}.
\end{align*}
Appending $n-k$ zeros to the principal curvatures preserves
$\sigma_m$ for $0\le m\le k$. Hence $\widetilde u_{\sigma_j}$ is
$k$-admissible and satisfies $G_{k,l}[\widetilde u_{\sigma_j}]
=\widetilde g$ in $B_\varepsilon^n$. The Newton--Maclaurin inequality
$(\sigma_k/\sigma_l)^{1/(k-l)}\le C(n,k,l)\sigma_k^{1/k}$
on $\Gamma_k$ shows that the quotient extends continuously
by zero outside $\Gamma_k$. Since the quotient is nondecreasing
in each variable on $\Gamma_k$ and
$\Gamma_k+[0,\infty)^n\subset\Gamma_k$, the extension remains
nondecreasing in each variable and thus defines a degenerate
elliptic operator. Because $\widetilde u_{\sigma_j}\to
\widetilde u_0$ uniformly, viscosity stability
\cite[Proposition~I.3]{Lions1983} gives
$G_{k,l}[\widetilde u_0]=\widetilde g$ in $B_\varepsilon^n$
in the viscosity sense. Finally, $\widetilde u_0$ is convex and
Lipschitz on $\overline{B_\varepsilon^n}$, and its restriction
$\widetilde u_0(x,0)=u_0(x)$ contradicts
$C^{1,\beta}(B_{\varepsilon/2}^n)$ regularity for every
$\beta>1-2/(k-l)$.

\hfill$\square$

\section*{Acknowledgements}
This work forms part of the author's ongoing doctoral research.
The author thanks Professor Siyuan Lu for helpful discussions
and for drawing his attention to related work. He is also grateful
to Dr. Bin Wang for valuable suggestions and to Professor Fei Han
and Yan Ma for helpful correspondence.

\section*{Use of artificial intelligence}

During the preparation of this manuscript, the author used GPT-6 Astra to assist with bibliographic searches, some mathematical calculations, the organization of proofs, and language editing. All AI-assisted suggestions were critically reviewed
by the author, who takes full responsibility for
the mathematical correctness and content of this paper.

\bibliography{main}

@misc{Wu2026,
  author = {Wu, Ke},
  title = {{Interior Hessian Estimates for $k$-Convex Solutions of Hessian Quotient Equations with $2k>n$}},
  year = {2026},
  note = {\href{https://arxiv.org/pdf/2609.23006}{arXiv:2609.23006}},
}

@misc{QiuYan2026,
  author = {Qiu, Guohuan and Yan, Jin},
  title = {{Interior Estimates for Hessian Quotient Equations}},
  year = {2026},
  note = {\href{https://arxiv.org/pdf/2609.24549}{arXiv:2609.24549}},
}

@misc{Fung2026,
  author = {Cheuk Yan Fung},
  title = {{Doubling Argument of the Hessian Estimate for the Hessian Quotient Equations}},
  year = {2026},
  note = {\href{https://arxiv.org/pdf/2607.21982}{arXiv:2607.21982}},
}

@article{ChenYuanWarren2009,
  author = {Chen, Jingyi and Warren, Micah and Yuan, Yu},
  title = {{A Priori Estimate for Convex Solutions to Special Lagrangian Equations and Its Application}},
  journal = {Communications on Pure and Applied Mathematics},
  year = {2009},
  volume = {62},
  number = {4},
  pages = {583--595},
  publisher = {Wiley Online Library},
  note = {\href{https://doi.org/10.1002/cpa.20261}{doi:\,\nolinkurl{10.1002/cpa.20261}}},
  mrnumber = {2492708},
}

@article{Qiu2024hessian,
  author = {Qiu, Guohuan},
  title = {{Interior Hessian Estimates for $\sigma_2$ Equations in Dimension Three}},
  journal = {Frontiers of Mathematics},
  year = {2024},
  volume = {19},
  number = {4},
  pages = {577--598},
  note = {\href{https://doi.org/10.1007/s11464-023-0156-0}{doi:\,\nolinkurl{10.1007/s11464-023-0156-0}}},
  mrnumber = {4768566},
}

@article{Fan2026D4,
  author = {Fan, Zhenyu},
  title = {Hessian estimates for the sigma-2 equation with variable right-hand side terms in dimension 4},
  journal = {Advances in Mathematics},
  year = {2026},
  volume = {494},
  publisher = {Elsevier BV},
  pages = {110953},
  note = {\href{https://doi.org/10.1016/j.aim.2026.110953}{doi:\,\nolinkurl{10.1016/j.aim.2026.110953}}},
  mrnumber = {5059337},
}

@misc{LiWu2026Sigma2,
  author = {Li, Zhisu and Wu, Ke},
  title = {Interior {Hessian} estimates for the quadratic {Hessian} equation},
  year = {2026},
  note = {\href{https://arxiv.org/pdf/2608.23233}{arXiv:2608.23233}},
}

@misc{DongZhang2026,
  author = {Dong, Weisong and Zhang, Ruijia},
  title = {Interior {Hessian} estimates for {Hessian} quotient equations},
  year = {2026},
  note = {\href{https://arxiv.org/pdf/2608.19087}{arXiv:2608.19087}},
}

@misc{LiWu2026Quotient,
  author = {Li, Zhisu and Wu, Ke},
  title = {A concavity inequality and interior {${C}^2$} estimate for {Hessian} quotient equations},
  year = {2026},
  note = {\href{https://arxiv.org/pdf/2608.17405}{arXiv:2608.17405}},
}

@article{Qiu2024curvature,
  author = {Qiu, Guohuan},
  title = {Interior curvature estimates for hypersurfaces of prescribing scalar curvature in dimension three},
  journal = {Amer. J. Math.},
  volume = {146},
  number = {3},
  pages = {579--605},
  year = {2024},
  note = {\href{https://doi.org/10.1353/ajm.2024.a928319}{doi:\,\nolinkurl{10.1353/ajm.2024.a928319}}},
  mrnumber = {4752675},
}

@article{Heinz1959,
  author = {Erhard Heinz},
  title = {{On Elliptic Monge--Amp\`ere Equations and Weyl's Embedding Problem}},
  journal = {Journal D’Analyse Mathematique},
  year = {1959},
  volume = {7},
  number = {1},
  pages = {1--52},
  publisher = {Springer},
  note = {\href{https://doi.org/10.1007/BF02787679}{doi:\,\nolinkurl{10.1007/BF02787679}}},
  mrnumber = {0111943},
}

@book{Pogorelov1978,
  title = {{The Minkowski Multidimensional Problem}},
  publisher = {Halsted Press [John Wiley \& Sons], New York-Toronto-London},
  year = {1978},
  author = {Aleksei Vasil'evich Pogorelov},
  series = {Scripta Series in Mathematics},
  address = {Washington D.C.},
  note = {Translated from the Russian by Vladimir Oliker; introduction by Louis Nirenberg; 106 pp},
  mrnumber = {0478079},
}

@article{CaffarelliNirenbergSpruck1988,
  author = {Luis A. Caffarelli and Louis Nirenberg and Joel Spruck},
  title = {{Nonlinear Second-Order Elliptic Equations. V. The Dirichlet Problem for Weingarten Hypersurfaces}},
  journal = {Communications on Pure and Applied Mathematics},
  volume = {41},
  number = {1},
  pages = {47--70},
  year = {1988},
  note = {\href{https://doi.org/10.1002/cpa.3160410105}{doi:\,\nolinkurl{10.1002/cpa.3160410105}}},
  mrnumber = {917124},
}

@article{Urbas1990,
  author = {John I. E. Urbas},
  title = {{On the Existence of Nonclassical Solutions for Two Classes of Fully Nonlinear Elliptic Equations}},
  journal = {Indiana University mathematics journal},
  year = {1990},
  pages = {355--382},
  publisher = {JSTOR},
  note = {\href{https://doi.org/10.1512/iumj.1990.39.39020}{doi:\,\nolinkurl{10.1512/iumj.1990.39.39020}}},
  mrnumber = {1089043},
}

@article{Trudinger1990,
  author = {Neil S. Trudinger},
  title = {{The Dirichlet Problem for the Prescribed Curvature Equations}},
  journal = {Archive for Rational Mechanics and Analysis},
  year = {1990},
  volume = {111},
  number = {2},
  pages = {153--179},
  publisher = {Springer},
  note = {\href{https://doi.org/10.1007/BF00375406}{doi:\,\nolinkurl{10.1007/BF00375406}}},
  mrnumber = {1057653},
}

@book{CaffarelliCabre1995,
  author = {Luis A. Caffarelli and Xavier Cabr{\'e}},
  title = {{Fully Nonlinear Elliptic Equations}},
  series = {American Mathematical Society Colloquium Publications},
  volume = {43},
  publisher = {American Mathematical Society},
  address = {Providence, RI},
  year = {1995},
  note = {\href{https://doi.org/10.1090/coll/043}{doi:\,\nolinkurl{10.1090/coll/043}}},
  mrnumber = {1351007},
}

@article{ShengUrbasWang2004,
  author = {Weimin Sheng and John I. E. Urbas and Xu-Jia Wang},
  title = {{Interior Curvature Bounds for a Class of Curvature Equations}},
  journal = {Duke Mathematical Journal},
  year = {2004},
  volume = {123},
  number = {2},
  pages = {235 -- 264},
  publisher = {Duke University Press},
  note = {\href{https://doi.org/10.1215/S0012-7094-04-12321-8}{doi:\,\nolinkurl{10.1215/S0012-7094-04-12321-8}}},
}

@article{WarrenYuan2009,
  author = {Micah Warren and Yu Yuan},
  journal = {Communications on Pure and Applied Mathematics},
  title = {{Hessian Estimates for the {$\sigma_2$} Equation in Dimension 3}},
  year = {2009},
  number = {3},
  pages = {305--321},
  volume = {62},
  note = {\href{https://doi.org/10.1002/cpa.20251}{doi:\,\nolinkurl{10.1002/cpa.20251}}},
  mrnumber = {2487850},
}

@article{WangYuan2014,
  author = {Dake Wang and Yu Yuan},
  title = {{Hessian Estimates for Special Lagrangian Equations with Critical and Supercritical Phases in General Dimensions}},
  journal = {American Journal of Mathematics},
  year = {2014},
  volume = {136},
  number = {2},
  pages = {481--499},
  publisher = {Johns Hopkins University Press},
  note = {\href{https://doi.org/10.1353/ajm.2014.0009}{doi:\,\nolinkurl{10.1353/ajm.2014.0009}}},
}

@article{GuanRenWang2015,
  author = {Pengfei Guan and Changyu Ren and Zhizhang Wang},
  title = {{Global {${C}^2$}-Estimates for Convex Solutions of Curvature Equations}},
  journal = {Communications on Pure and Applied Mathematics},
  year = {2015},
  volume = {68},
  number = {8},
  pages = {1287-1325},
  note = {\href{https://doi.org/10.1002/cpa.21528}{doi:\,\nolinkurl{10.1002/cpa.21528}}},
  mrnumber = {3366747},
}

@article{McGonagleSongYuan2019,
  author = {Matt McGonagle and Chong Song and Yu Yuan},
  title = {{Hessian Estimates for Convex Solutions to Quadratic Hessian Equation}},
  journal = {Annales de l'Institut Henri Poincar{\'e} C, Analyse non lin{\'e}aire},
  year = {2019},
  volume = {36},
  number = {2},
  pages = {451--454},
  note = {\href{https://doi.org/10.1016/j.anihpc.2018.07.001}{doi:\,\nolinkurl{10.1016/j.anihpc.2018.07.001}}},
  mrnumber = {3913193},
}

@article{ShankarYuan2020,
  author = {Ravi Shankar and Yu Yuan},
  title = {{Hessian Estimate for Semiconvex Solutions to the {$\sigma_2$} Equation}},
  journal = {Calculus of Variations and Partial Differential Equations},
  year = {2020},
  volume = {59},
  number = {1},
  pages = {Paper No. 30, 12 pp.},
  publisher = {Springer},
  note = {\href{https://doi.org/10.1007/s00526-019-1690-1}{doi:\,\nolinkurl{10.1007/s00526-019-1690-1}}},
  mrnumber = {4054864},
}

@article{Mooney2021,
  author = {Connor Mooney},
  title = {Strict {$2$}-convexity of convex solutions to the quadratic {Hessian} equation},
  journal = {Proceedings of the American Mathematical Society},
  year = {2021},
  volume = {149},
  number = {6},
  pages = {2473--2477},
  note = {\href{https://doi.org/10.1090/proc/15454}{doi:\,\nolinkurl{10.1090/proc/15454}}},
  mrnumber = {4246798},
}

@misc{Lu2023,
  author = {Siyuan Lu},
  title = {{Interior {${C}^2$} Estimate for Hessian Quotient Equation in Dimension Three}},
  year = {2023},
  note = {\href{https://arxiv.org/pdf/2311.05835}{arXiv:2311.05835}},
}

@article{Zhou2024,
  author = {Xingchen Zhou},
  title = {{Notes on Generalized Special Lagrangian Equation}},
  journal = {Calculus of Variations and Partial Differential Equations},
  year = {2024},
  volume = {63},
  number = {8},
  pages = {197},
  publisher = {Springer},
  note = {\href{https://doi.org/10.1007/s00526-024-02801-w}{doi:\,\nolinkurl{10.1007/s00526-024-02801-w}}},
}

@article{ShankarYuan2025,
  author = {Ravi Shankar and Yu Yuan},
  journal = {Annals of Mathematics},
  title = {{Hessian Estimates for the {$\sigma_2$} Equation in Dimension Four}},
  year = {2025},
  number = {2},
  pages = {489 -- 513},
  volume = {201},
  publisher = {Department of Mathematics of Princeton University},
  note = {\href{https://doi.org/10.4007/annals.2025.201.2.4}{doi:\,\nolinkurl{10.4007/annals.2025.201.2.4}}},
  mrnumber = {4880431},
}

@article{Lu2025,
  author = {Siyuan Lu},
  title = {{Interior {${C}^2$} Estimate for Hessian Quotient Equation in General Dimension}},
  journal = {Ann. PDE},
  year = {2025},
  volume = {11},
  pages = {17},
  note = {\href{https://doi.org/10.1007/s40818-025-00215-1}{doi:\,\nolinkurl{10.1007/s40818-025-00215-1}}},
}

@misc{FanShankar2026,
  author = {Zhenyu Fan and Ravi Shankar},
  title = {{Interior Estimates and Regularity for the Scalar Curvature Equation in Dimension 4}},
  year = {2026},
  note = {\href{https://arxiv.org/pdf/2608.22748}{arXiv:2608.22748}},
}

@article{Liu2026,
  author = {Jianxiang Liu},
  title = {{Interior Curvature Estimate for Curvature Quotient Equations on Convex Hypersurfaces}},
  journal = {Communications in Analysis and Geometry},
  year = {2026},
  note = {To appear. \href{https://arxiv.org/pdf/2505.00360}{https://arxiv.org/pdf/2505.00360}},
}

@article{guanqiu2019,
  author = {Pengfei Guan and Guohuan Qiu},
  title = {{Interior ${C}^2$ regularity of convex solutions to prescribing scalar curvature equations}},
  year = {2019},
  journal = {Duke Mathematical Journal},
  volume = {168},
  number = {9},
  pages = {1641--1663},
  note = {\href{https://doi.org/10.1215/00127094-2019-0001}{doi:\,\nolinkurl{10.1215/00127094-2019-0001}}},
  mrnumber = {3961212},
}

@article{chenjianz2026,
  author = {Ruosi Chen and Huaiyu Jian and Xingchen Zhou},
  title = {{An integral approach to prescribing scalar curvature equations}},
  year = {2026},
  journal = {Mathematische Annalen},
  volume = {394},
  pages = {Article 72},
  note = {\href{https://doi.org/10.1007/s00208-026-03312-z}{doi:\,\nolinkurl{10.1007/s00208-026-03312-z}}},
}

@article{shankaryuan2021,
  author = {Ravi Shankar and Yu Yuan},
  title = {{Regularity for almost convex viscosity solutions of the sigma-2 equation}},
  year = {2021},
  journal = {Journal of Mathematical Study},
  volume = {54},
  number = {2},
  pages = {164--170},
  note = {\href{https://doi.org/10.4208/jms.v54n2.21.03}{doi:\,\nolinkurl{10.4208/jms.v54n2.21.03}}},
}

@article{mooney2025,
  author = {Connor Mooney},
  title = {{Remarks on the quadratic Hessian equation}},
  year = {2025},
  journal = {J. Math. Study},
  note = {To appear. \href{https://arxiv.org/abs/2505.14586}{https://arxiv.org/abs/2505.14586}},
}

@misc{chenjiantuzhou2026,
  author = {Ruosi Chen and Huaiyu Jian and Xushan Tu and Xingchen Zhou},
  title = {{Regularity for convex viscosity solutions of $\sigma_2$ equation}},
  year = {2026},
  note = {\href{https://arxiv.org/pdf/2605.30823}{arXiv:2605.30823}},
}

@misc{zhouzhu2026,
  author = {Xingchen Zhou and Ruixuan Zhu},
  title = {{Interior ${C}^{2,\alpha}$ regularity for convex solutions of the $2$-Hessian equation}},
  year = {2026},
  note = {\href{https://arxiv.org/pdf/2608.24604}{arXiv:2608.24604}},
}

@misc{chenzhouzhu2026,
  author = {Ruosi Chen and Xingchen Zhou and Ruixuan Zhu},
  title = {{Interior ${C}^{2,\alpha}$ regularity for the quadratic Hessian equation}},
  year = {2026},
  note = {\href{https://arxiv.org/pdf/2608.29484}{arXiv:2608.29484}},
}

@misc{qiuyan2026curvature,
  author = {Guohuan Qiu and Jin Yan},
  title = {{Interior curvature estimates for the graphical scalar curvature equation in all dimensions}},
  year = {2026},
  note = {\href{https://arxiv.org/pdf/2609.02581}{arXiv:2609.02581}},
}

@article{warrenyuan2010,
  author = {Warren, Micah and Yuan, Yu},
  title = {Hessian and gradient estimates for three dimensional special {Lagrangian} equations with large phase},
  journal = {Amer. J. Math.},
  volume = {132},
  number = {3},
  pages = {751--770},
  year = {2010},
  note = {\href{https://doi.org/10.1353/ajm.0.0115}{doi:\,\nolinkurl{10.1353/ajm.0.0115}}},
  mrnumber = {2666907},
}

@article{lutsai2026note,
  author = {Siyuan Lu and Yi-Lin Tsai},
  title = {{A note on interior ${C}^2$ estimate for general Hessian quotient equation}},
  year = {2026},
  journal = {Communications on Pure and Applied Analysis},
  volume = {33},
  pages = {88--100},
  note = {\href{https://doi.org/10.3934/cpaa.2026045}{doi:\,\nolinkurl{10.3934/cpaa.2026045}}},
}

@misc{jiaosui2026,
  author = {Heming Jiao and Zhenan Sui},
  title = {{Interior Hessian estimates for Hessian quotient equations in dimension three}},
  year = {2026},
  note = {\href{https://arxiv.org/pdf/2602.14064}{arXiv:2602.14064}},
}

@misc{lusroka2026,
  author = {Siyuan Lu and Marcin Sroka},
  title = {{On Liouville's theorem for the Hessian quotient equation $\sigma_2/\sigma_1$}},
  year = {2026},
  note = {\href{https://arxiv.org/pdf/2602.14946}{arXiv:2602.14946}},
}

@misc{qiuzhou2024,
  author = {Guohuan Qiu and Xingchen Zhou},
  title = {{A priori interior estimates for special Lagrangian curvature equations}},
  year = {2024},
  note = {\href{https://arxiv.org/pdf/2407.15159}{arXiv:2407.15159}},
}

@misc{hanma2026,
  author = {Fei Han and Yan Ma},
  title = {{Interior curvature estimates for graphical curvature quotient equations}},
  year = {2026},
  note = {\href{https://arxiv.org/pdf/2609.39276}{arXiv:2609.39276}},
}

@article{urbas2000curvature,
  author = {Urbas, John I. E.},
  title = {An interior curvature bound for hypersurfaces of prescribed $k$-th mean curvature},
  journal = {J. Reine Angew. Math.},
  volume = {519},
  pages = {41--57},
  year = {2000},
  note = {\href{https://doi.org/10.1515/crll.2000.016}{doi:\,\nolinkurl{10.1515/crll.2000.016}}},
}

@misc{jiliang2026,
  author = {Ji, Ke and Liang, Lichun},
  title = {Interior {Hessian} estimates for semi-convex solutions of the $\sigma_2/\sigma_1$ equation with {Lipschitz} right-hand sides},
  year = {2026},
  note = {\href{https://arxiv.org/pdf/2608.24530}{arXiv:2608.24530}},
}

@article{Lions1983,
  title = {Optimal control of diffusion processes and hamilton–jacobi–bellman equations part 2 : viscosity solutions and uniqueness},
  volume = {8},
  issn = {1532-4133},
  url = {http://dx.doi.org/10.1080/03605308308820301},
  doi = {10.1080/03605308308820301},
  number = {11},
  journal = {Communications in Partial Differential Equations},
  publisher = {Informa UK Limited},
  author = {Lions, P. L.},
  year = {1983},
  month = Jan,
  pages = {1229–1276},
  mrnumber = {709162},
}

@misc{liwu2026generalquotient,
  author = {Zhisu Li and Ke Wu},
  title = {{Interior Hessian estimates for $k$-convex solutions of Hessian quotient equations}},
  year = {2026},
  note = {\href{https://arxiv.org/pdf/2610.18252}{arXiv:2610.18252}},
}

@article{Wang2026CurvatureQuotient,
  author = {Bin Wang},
  title = {{Interior curvature estimates for convex solutions of a prescribed curvature quotient equation}},
  year = {2026},
  note = {To appear. \href{https://arxiv.org/abs/2608.20188}{arXiv:2608.20188}},
}
\bibliographystyle{amsplain}
\end{document}